\documentclass[12pt]{article}

\usepackage{amsmath,amsthm}
\usepackage{amssymb}
\usepackage{amsfonts}
\usepackage{amscd}

\usepackage{color}

\usepackage[normalem]{ulem}

\usepackage{cancel}

\newtheorem{thm}{Theorem}[section]
\newtheorem{theorem}[thm]{Theorem}

\newtheorem{corollary}[thm]{Corollary}
\newtheorem{lemma}[thm]{Lemma}
\newtheorem{proposition}[thm]{Proposition}

\theoremstyle{definition}
\newtheorem{definition}[thm]{Definition}
\newtheorem{example}[thm]{Example}
\newtheorem{examples}[thm]{Examples}
\newtheorem{remark}[thm]{Remark}

\newtheorem{observation}[thm]{Observation}

\newcommand{\N} {\mathbb{N}}
\newcommand{\Z} {\mathbb{Z}}

\newcommand{\C} {\mathbb{C}}
\newcommand{\KK} {\mathbb{K}}
\newcommand{\bk} {\Bbbk}

\newcommand{\bt} {\mathbf{t}}
\newcommand{\bu} {\mathbf{u}}
\newcommand{\bG} {\mathbf{G}}

\newcommand{\bA} {\mathbf{A}}

\newcommand{\fg} {\mathfrak{g}}
\newcommand{\fh} {\mathfrak{h}}
\newcommand{\fH} {\mathcal{H}}

\newcommand{\kzg} {K_{\Z}(\fg)} 
\newcommand{\uzg} {U_{\Z}(\fg)} 

\newcommand{\salg} {\mathrm{(salg)}}

\newcommand{\sets}{{(\mathrm{sets})}}
\newcommand{\grps} {\mathrm{(groups)}}
\newcommand{\spec}{{\hbox{\sl Spec}\,}}
\newcommand{\uspec}{\underline{\hbox{\sl Spec}\,}}

\newcommand{\Hom}{\mathrm{Hom}}

\newcommand{\lra} {\longrightarrow}

\newcommand{\End}{\mathrm{End}}

\newcommand{\rGL}{\mathrm{GL}}

\newcommand{\Lie}{\mathrm{Lie}}
\newcommand{\str}{\mathrm{str}}
\newcommand{\tr}{\mathrm{tr}}

\newcommand{\rosp}{\mathfrak{osp}}
\newcommand{\rsl}{\mathfrak{sl}}
\newcommand{\rgl}{\mathfrak{gl}}
\newcommand{\rsp}{\mathfrak{sp}}
\newcommand{\rso}{\mathfrak{so}}

\newcommand{\cO}{\mathcal{O}}
\newcommand{\cC}{\mathcal{C}}

\newcommand{\wDelta}{{\Delta}}

\newcommand{\al}{{\alpha}}
\newcommand{\be}{{\beta}}
\newcommand{\ga}{{\gamma}}

\newcommand{\ep}{{\epsilon}}

\begin{document}

\bigskip

\centerline{\Large \bf ON THE CONSTRUCTION OF}

\bigskip

\centerline{\Large \bf CHEVALLEY SUPERGROUPS}

\bigskip

\centerline{R. Fioresi$^\flat$, F. Gavarini$^\#$}

\bigskip

\centerline{\it $^\flat$ Dipartimento di Matematica, Universit\`a
di Bologna } \centerline{\it Piazza di Porta San Donato, 5 ---
I-40127 Bologna, Italy} \centerline{{\footnotesize e-mail:
fioresi@dm.unibo.it}}

\bigskip

\centerline{\it $^\#$ Dipartimento di Matematica, Universit\`a di
Roma ``Tor Vergata'' } \centerline{\it via della ricerca
scientifica 1  --- I-00133 Roma, Italy}

\centerline{{\footnotesize e-mail: gavarini@mat.uniroma2.it}}

\begin{abstract}
We give a description
of the construction of Chevalley supergroups, providing
some explanatory examples. We avoid the discussion of the $A(1,1)$, $P(3)$ and
$Q(n)$
cases, for which our construction holds, but the exposigetion becomes more
complicated.
We shall not in general provide complete proofs for our statements,
instead we will make an effort to convey the key ideas underlying our
construction. A fully detailed account of our work is scheduled to
appear in \cite{fg}.
\end{abstract}

\section{Introduction}

The notion of Chevalley group, introduced by Chevalley in 1955,
provided a unified combinatorial construction of all
simple algebraic groups over a generic field $k$. The consequences
of Chevalley's work were many and have had tremendous impact
in the following decades. 
His construction was motivated by issues linked to the problem
of the classification of semisimple algebraic groups: he provided
an existence theorem for such groups, essentially exhibiting an
example of simple group for each of the predicted possibility.
In the course of this discussion, he discovered new examples of 
finite simple groups, which
had escaped to the group theorists up to then. 
Later on, in the framework of a modern treatment of algebraic geometry,
his work was instrumental to show that all simple algebraic
groups are algebraic schemes over $\Z$ and to study arithmetic
questions over arbitrary fields.

\medskip

We may say that we have similar motivations: we want a unified approach
to describe all algebraic supergroups, which have Lie superalgebras of
classical type and we also want to give new examples of supergroups,
over arbitrary fields. For instance, our discussion enables us to 
provide an explicit
construction of algebraic supergroups associated
with the exceptional and the strange Lie superalgebras. To our knowledge
these supergroups have not been examined before, though an approach  
in the differential setting can be very well carried through
via the language of {\sl super Harish-Chandra pairs}. 
In such approach a supergroup is  understood as a pair ($G_0$, $\fg$), 
consisting of an ordinary group $G_0$ and a super Lie algebra $\fg$, 
with even part $\fg_0=\Lie(G_0)$, together with some natural
compatibility conditions involving the adjoint action of the
group $G_0$ on $\fg$.
It is clear that in positive characteristic this method shows severe
limitations.

\medskip

In the present work we outline the construction of the
Chevalley supergroups associated with Lie superalgebras of classical type.
We shall not present complete proofs for our statements, 
they will appear in \cite{fg}, however we shall
concentrate on the key ideas and examples that will help to understand our
construction.
 
In our statements, we shall leave out the 
strange Lie superalgebra $Q(n)$ and some low dimensional
cases, which can be treated very well with the same method, 
with minor modifications,
but present extra difficulties that make our construction 
and notation opaque.

\medskip

Essentially, we are going to follow Chevalley's recipe and push it
as far as we can, before resorting to more sophisticated algebraic
geometry tecniques,
when the supergeometric nature of our objects forces us to do so.

\medskip

We start with a complex 
Lie superalgebra of classical type $\fg$, together with a 
fixed Cartan subalgebra $\fh$,
and we define the {\sl Chevalley basis} of $\fg$. This is an homogeneous 
basis of $\fg$, as super vector space, whose elements have the brackets
expressed as a linear combination of the basis elements 
with just {\sl integral}
coefficients. Consequently they
give us an integral form of $\fg$, that we call $\fg_\Z$
the {\sl Chevalley Lie superalgebra} associated with $\fg$ and $\fh$.
Such integral form gives raise to the Kostant integral form $\kzg$
of the universal enveloping superalgebra $U(\fg)$ 
of $\fg$. $\kzg$ is free over $\Z$ with
basis given by the ordered monomials in the divided powers of the
root vectors and the binomial coefficients
in the generators of $\fh$ in the Chevalley basis: 
$X^m/m!$, $\left( \begin{array}{c} H_i \\ n 
\end{array}\right)$, $\al \in \Delta$ (root system) and $m,n \in \N$.

\medskip

Next, we look at a faithful rational representation of $\fg$ in a  
finite dimensional complex vector space $V$. Inside $V$ we can
find an {\sl integral lattice} $M$ which is invariant under
the action of $\kzg$ and its stabilizer $\fg_V$ in $\fg$ defines an
integral form of $\fg$. In complete
analogy with Chevalley, for an arbitrary field $k$, we can give the following
key definitions:
$$
V_k:=k \otimes_\Z M, \qquad \fg_k:=k \otimes_\Z \fg_V, \qquad
U_k:=k \otimes_\Z \kzg. 
$$
We could even take $k$ to be a commutative ring, however for the
scope of the present work and to stress the analogy with Chevalley's
construction, we prefer the restrictive hypothesis of $k$ to be a field.

\medskip

This is the point where our construction departs dramatically from
Chevalley's one. In fact, starting from the faithful representation
$V_k$ of $\fg_k$, 
Chevalley defines the Chevalley group $G_V$ as generated by the
exponentials 
$exp(tX_\al):=1+tX_\al+(t^2/2) X_\al^2+\dots$, for $t \in k$ and 
$X_\al$ the root vector corresponding
to the root $\alpha$ in the Chevalley basis. 
Such an expression makes sense
since the $X_\al$'s act as nilpotent elements. If we were to repeat without
changes this construction in the super setting, we shall 
find only ordinary groups over $k$
associated with the Lie algebra $\fg_0$, 
the even part of $\fg$.
This is because over a field, we cannot see any supergeometric
behaviour; the only thing we can recapture is the underlying 
classical object. For this reason, we need to go beyond Chevalley's
construction and build our supergroups as {\sl functors}.

\medskip

We define $\bG$ the {\sl Chevalley supergroup} associated
with $\fg$ and the faithful representation $V$, as the
functor $\bG:\salg \lra \sets$, with $\bG(A)$ the subgroup
of $\rGL(A \otimes V_k)$ generated by $\bG_0(A)$ and the elements 
$1+\theta_\beta X_\beta$, for $\beta \in \Delta_1$. 
In other words we have:
$$
\bG_V(A)=\langle \bG_0(A), 1+\theta_\beta X_\beta \rangle \subset
\rGL(A \otimes V_k), \qquad A \in \salg, \qquad \theta_\beta \in A_1
$$
where $\salg$ and $\sets$ are the categories of commutative
superalgebras and sets respectively and (as always) we
use $X_\beta$ to denote also the
image of the root vector $X_\beta$ 
in the chosen faithful representation $V_k$. $\bG_0$ is the
functor of points of the (reductive) algebraic supergroup associated to
$\fg_0$ and the representation $V_k$.

\medskip

This is a somehow  natural generalization of what Chevalley does in his
original construction: he provides
the $k$-points of the algebraic group scheme constructed starting from
a complex semisimple Lie algebra and a faithful representation,
for all the fields $k$,
while we give the $A$-points of the supergroup scheme for any commutative
$k$-superalgebra $A$.

\medskip

Once this definition is properly established, we need to show
that $\bG$ is the functor of points of an algebraic supergroup,
in other words, that it is representable. This is the price to
pay when we employ the language of the functor of points: it is
much easier to define geometric objects, however we need to prove
representability in order to speak properly of supergroup schemes.
As customary, we use the same letter to denote both the superscheme
and its functor of points.
\medskip
 
We shall obtain the representability of $\bG$ by showing that
$$
\bG \cong \bG_0 \times \bA^{0|N}
$$ 
where $\bA^{0|N}$ is the functor of points of an affine superspace
of dimension $0|N$. Once this isomorphism is established the
representability follows at once, since both $\bG_0$ and
$\bA^{0|N}$ are representable, i.e. they are the functors of points of
superschemes, hence their product is.

\medskip

The next question we examine is how much our construction
depends on the chosen representation.
In complete analogy to Chevalley approach, we show that
if we have two representations $V$ and $V'$, with weight lattices
$L_V \subset L_{V'}$, then there is a surjective morphism $\bG_{V'} \lra 
\bG_{V}$, with kernel in the center of $\bG_{V'}$. This implies right
away that our construction depends only on the weight lattice of
the chosen representation $V$ and in particular it shows that it
is independent from the choice of
the lattice $M$ inside $V$.
\medskip

This paper is organized as follows.

\medskip

In section \ref{supergeometry} we review quickly some facts of
algebraic supergeometry and the theory of  Lie superalgebras.

\medskip

In sections \ref{chevalleybasis} and \ref{chesupergroups} we go to the heart of the construction
of Chevalley's supergroups going through all
the steps detailed above.

\medskip

Finally in section \ref{examples} we provide some insight into our construction
with some examples and observations.

\medskip

We wish to thank the UCLA Mathematics Department, for hosting the
conference. 
We thank I. Dimitrov, V. Serganova and A. Schwarz for
helpful comments. We finally thank our referee for helping us to
improve the readibility of our manuscript.

\section{Supergeometry in the algebraic setting} \label{supergeometry}

Let $k$ be the ground field.

\medskip

A \textit{super vector space} $V$ is a vector space with $\Z_2$ grading:
$V=V_0 \oplus V_1$, the elements in $V_0$ are called {\it even}
and the elements in $V_1$ are called {\it odd}. 
Hence we have a function $p$ 
called the \textit{parity} defined only on homogeneous elements.
A \textit{superalgebra} $A$ is a super vector space with multiplication
preserving parity;
$A$ is \textit{commutative} if $xy=(-1)^{p(x)p(y)}yx$ for all
$x$, $y$ homogeneous elements in $A$. All superalgebras are assumed to be
commutative unless otherwise specified and their category is
denoted with $\salg$.

\begin{definition}
A  {\it superspace}  $ \, S = \big( |S|, \cO_S \big) \, $  is a 
topological space  $ |S| $  endowed with a sheaf of commutative 
superalgebras  $ \cO_S $  such that the stalk  $ \cO_{S,x} $  is 
a local superalgebra for all  $ \, x \in |S| \, $.

\medskip

A  {\it morphism}  $ \; \phi: S \lra T \; $  of superspaces 
consists of a pair $ \; \phi = \big( |\phi|, \phi^* \big) \, $,  
where  $ \; \phi : |S| \lra |T| \; $  is a morphism of topological 
spaces and  $ \; \phi^* : \cO_T \lra \phi_* \cO_S \; $  is a sheaf 
morphism such that  $ \; \phi_x^* \big( {\mathbf{m}}_{|\phi|(x)} \big)
= {\mathbf{m}}_x \; $  where  $ {\mathbf{m}}_{|\phi|(x)} $ and  $
{\mathbf{m}}_{x} $  are the maximal ideals in the stalks $ \cO_{T,
\, |\phi|(x)} $  and  $ \cO_{S,x} $  respectively and $\phi_x^*$
is the morphism induced by $\phi^*$ on the stalks
and $ \phi_* \cO_S $  is the sheaf on  $ |T| $  
defined as  $ \, \phi_* \cO_S(V) := \cO_S(\phi^{-1}(V)) \, $.
\end{definition}

The next example of superspace 
turns out to be extremely important,
as $\uspec A$, for a commutative superalgebra $A$,
is the local model for superschemes, very much
in the same way as $\uspec A_0$ is the local model for
ordinary schemes for $A_0$ a commutative algebra. 

\begin{example}
Let  $ \, A \in \salg \, $  and let  $ \cO_{A_0} $  
be the structural sheaf of the ordinary scheme  
$ \, \uspec(A_0) = \big( \spec(A_0), \cO_{A_0} \big) \, $,  
where  $ \spec(A_0) $  denotes the prime spectrum of the commutative ring  
$ A_0 \; $.  Now  $ A $  is a module over  $ A_0 \, $,  so we have a 
sheaf  $ \cO_A $  of  $ \cO_{A_0} $--modules  over  $ \spec(A_0) $  
with stalk  $ A_p \, $,  the  $ p $--localization of the  $ A_0 $--module  
$ A \, $,  at the prime  $ \, p \in \spec(A_0) \, $.

\medskip

$ \; \uspec(A) := \big( \spec(A_0), \cO_A \big) \; $ is
a superspace, as one can readily check.

\medskip

Given $f:A \lra B$ a superalgebra morphism, one can define
$\uspec f: \uspec B \lra \uspec A$ in a natural way, very similarly to
the ordinary setting, thus making $\uspec$ a
functor $\uspec: \salg \lra \sets$, where $\salg$ is the category of
superalgebras and $\sets$ the category of sets (see \cite{cf} ch. 5 
or \cite{eh} ch. 1, for more details).
\end{example}

\begin{definition}
Given a superspace $X$, we say it is an {\it affine superscheme}
if it is isomorphic to  $ \, \uspec(A) \, $  
for some commutative superalgebra  $ A \, $. We say that $X$ is a
{\it superscheme} if it is locally isomorphic to an affine superscheme.
\end{definition}

\begin{example} \label{affinesuperspace}
The  {\sl affine superspace\/}  $ \, \mathbb{A}_\bk^{p|q} \, $,  
also denoted  $ \, \bk^{p|q} \, $,  is defined  
as  
$$
\mathbb{A}_\bk^{p|q} := \bk[x_1 \dots x_p]
\otimes \wedge (\xi_1 \dots \xi_q) 
$$
where $\wedge (\xi_1 \dots \xi_q)$  is the exterior  
algebra generated by  the indeterminates $ \xi_1 $,  $ \dots $,  $ \xi_q \, $.
\end{example}

The formalism of the functor of points that we borrow from
algebraic geometry allows us to handle supergeometric objects
that would be otherwise very difficult to treat using just the
superschemes language. 

\begin{definition}
Let  $ X $  be a superscheme.  Its  {\it functor of points\/}
is the functor defined on the objects as
  $$  h_X : \salg \lra \sets \;\; ,
\qquad  h_X(A) \, := \, \Hom \big(\, \uspec(A) \, , X \big)  $$
and on the arrows as  $ \; h_X(f)(\phi) := \phi \circ \uspec (f) \; $.

\medskip

Since the category of affine superschemes is equivalent to the
category of commutative superalgebras (ref. \cite{cf, eh}) we have
that, when $X$ is affine, its functor of points is equivalently
defined as follows:
$$
h_X(A)=\Hom \big( \cO(X) \, , A \big), \qquad h_X(f)(\phi)=f \circ \phi
$$
where $\cO(X)$ is the superalgebra of global
sections of the structure sheaf on  $ X $.

\medskip

If  $ h_X $  is group valued, i.~e.~it is valued
in the category  $ \grps $  of groups, we say that  $ X $  is a  {\it
supergroup}. When  $ X $  is affine, this is equivalent to the
fact that $ \, \cO(X) \, $   is a
(commutative)  {\sl Hopf superalgebra}. More in general, we call  
{\it supergroup functor\/}  any functor  $ \; G : \salg \lra \grps \; $.
                                                
\medskip

   Any representable supergroup functor is the same as an affine supergroup.
   Following a customary abuse of notation, we shall then use the same letter 
to denote both the superscheme  $ X $  and its functor of points  $ h_X \, $.

\end{definition}

As always, Yoneda's lemma plays a crucial role,
allowing us to use natural transformations between the functors of points
of superschemes and the morphisms of  the superschemes 
themselves interchangeably.

\begin{proposition} ({\sl Yoneda's Lemma}) 
Let  $ \cC $  be a category,  
and  let $ R $,  $ S $ be two objects in  $ \cC \, $.  
Consider the two functors
$ \; h_R \, , \, h_S : \cC \lra \sets \; $  
defined on the objects by  $ \; h_R(A)
 := \Hom \big(\, R , A \big) \, $,  $ \; 
h_S(A) := \Hom \big(\, S , A \big) \; $  and on the
arrows by  $ \; h_R(f)(\phi) := f \circ \phi \, $,  $ \; 
h_S(f)(\psi) := f \circ \psi \; $.
                                                      
Then there exists a one-to-one correspondence between the natural 
transformations and the morphisms
$$ 
\, \big\{ h_R \lra h_S \big\} \, \longleftrightarrow \, \Hom(R,S) \,.
$$

\end{proposition}

This has an immediate corollary.

\begin{corollary}
Two affine superschemes are isomorphic if and only if
their functors of points are isomorphic.
\end{corollary}

   The next examples turn out to be very important in the sequel.

\medskip

\begin{examples} {\ }

 \vskip4pt

{\it (1)} \,  {\sl Super vector spaces as superschemes}.
Let  $ V $  be a super vector space.  
For any superalgebra  $ A $  we define  $ \; V(A) \, := 
\, {(A \otimes V)}_0 \, = \, A_0 \otimes V_0 \oplus A_1 \otimes V_1 \; $.  
This is a representable functor in the category of superalgebras, 
whose representing object is  $ {\hbox{Pol}}(V) \, $,  
the algebra of polynomial functions on  $ V \, $. 
Hence any super vector space can
be equivalently viewed as an affine superscheme.
If $V=k^{m|n}$, that is 
$ \, V_0 \cong k^p \, $  and  $ \, V_1 \cong k^q \, $,
$V$ is the functor of points of the affine superspace described
in \ref{affinesuperspace}.
\medskip

{\it (2)} \,  {\sl  $ \rGL(V) $  as an algebraic supergroup}.  
Let  $ V $  be a finite dimensional super vector space of dimension $p|q$.  
For any superalgebra  $ A \, $,  
let  $ \, \rGL(V)(A) := \rGL\big(V(A)\big) \, $  be the set of isomorphisms  
$ \; V(A) \lra V(A) \; $.  If we fix a homogeneous basis for  $ V \, $,  
we see that  $ \, V \cong k^{p|q} \; $.  
In this case, we also denote $ \, \rGL(V) \, $  with  $ \, \rGL(p|q) \, $.  
Now,  $ \rGL(p|q)(A) $  is the group of invertible matrices of size  
$ (p+q) $  with diagonal block entries in  $ A_0 $  and off-diagonal block 
entries in $ A_1 \, $.  It is well known that the functor  $ \rGL(V) $  
is representable; see (e.g.),  \cite{vsv}, Ch.~3,  for further details.

\end{examples}

We end our minireview of supergeometry by introducing the
concept of {\sl Lie superalgebra} and stating the Kac's classification
theorem for Lie superalgebras of classical type.

\medskip

We assume now  $\textit{char}(k) \neq 2,3$.

\begin{definition}
Let  $ \, \fg = \fg_0 \oplus \fg_1 \, $  be a super 
vector space.  
We say that  $ \fg $  is a Lie superalgebra, 
if we have a bracket  $ \; [\ ,\ ] : \fg \times \fg \lra \fg \; $  
which satisfies the following properties 
(as usual for all  $ x $,  $ y \in \fg \, $  homogeneous):

\begin{itemize}

\item {\sl Anti-symmetry:}  
$$
[x,y] \, + \, {(-1)}^{p(x) \, p(y)}[y,x] \; = \; 0 
$$

\item {\sl Jacobi identity:}
$$  
{(-\!1)}^{p(x) \, p(z)} \, [x,[y,z]] \, + \, {(-\!1)}^{p(y) \, p(x)} \, 
[y,[z,x]] \, + \, {(-\!1)}^{p(z) \, p(y)} \, [z,[x,y]] \; = \; 0.  
$$

\end{itemize}

\end{definition}

\medskip

The standard example is $\End(V)$   
the endomorphisms of the super vector space $V$, with $\End(V)_0$ the endomorphisms
preserving parity and $\End(V)_1$ the endomorphisms reversing parity.
The bracket is defined as:
$$
[X,Y]:=XY-(-1)^{|X||Y|}YX.
$$
If $ \, V := k^{p|q} = k^p \oplus k^q \, $,  
with  $ \, V_0 := k^p \, $  and  $ \, V_1 := k^q \, $,
we write  $ \; \End\big(k^{p|q}\big) := \, \End(V) \; $  
or  $ \; \rgl(p\,|q) := \End(V) \; $.
In this case $\End(V)_0$ consists of diagonal block matrices,
while $\End(V)_1$ consists of off diagonal block matrices all
with entries in $k$.

\medskip
                                                                     
In  $ \End(k^{p|q}) $  we can define the  {\it supertrace\/}  as follows:
$$  \str \begin{pmatrix}
              A  &  B  \\
              C  &  D
           \end{pmatrix}  := \;  \tr(A) - \tr(D)  \quad .  
$$

\medskip

There is an important class of Lie superalgebras, namely the
{\sl simple Lie superalgebras} that have been classified by Kac
(see \cite{ka}).

\begin{definition}
A non abelian Lie superalgebra $\fg$ is \textit{simple} if it has no nontrivial
homogeneous ideals. A Lie superalgebra $\fg$ is called of
\textit{classical type \/}  if it is simple and $ \fg_1 $  is 
completely reducible as a  $ \fg_0 $--module.  
Furthermore,  $ \fg $  is said to be  \textit{basic\/}  if, 
in addition, it admits a non-degenerate, invariant bilinear form.%
\end{definition}

\medskip

We now give a list of Lie superalgebras of classical type, sending the
reader to \cite{ka}, \cite{sc} for the details.

\begin{examples} 

\medskip

   {\it (1)} ---  $ \rsl(m|\,n) \, $.  
Define  $ \, \rsl(m|n) \, $  as the subset of  $ \rgl(m|n)  $  consisting of
all matrices with supertrace zero.  This is a Lie subalgebra of  
$ \rgl(m|n) \, $,  with  the following $ \Z_2 $--grading:
$$  
\rsl(m|\,n)_0 \, = \, \rsl(m) \oplus \rsl(n) \oplus \rgl(1)  \;\; ,  \qquad
\rsl(m|\,n)_1 \, = \, f_m \otimes f'_n \, \oplus \, f'_m \otimes f_n  
$$
where  $ f_r $  is the defining representation of  $ \rsl(r) $
and  $ f'_r $  is its dual (for any  $ r \, $).  When  $ \, m \neq
n \, $  this is a Lie superalgebra of {\sl classical type\/}.

\medskip

   {\it (2)} ---  $ \rosp(p\,|q) \, $.  
Let  $ \phi $  denote a nondegenerate consistent supersymmetric bilinear 
form in  $ \, V := k^{p|q} \, $.  
This means that  $ V_0 $  and  $ V_1 $  are mutually orthogonal and the 
restriction of  $ \phi $  to  $ V_0 $  is a symmetric and to  $ V_1 $  a 
skewsymmetric form (in particular,  $ \, q = 2\,n \, $  is even).  
We define in  $ \rgl(p\,|q) $  the subalgebra  
$ \; \rosp(p\,|q) := {\rosp(p ,|q)}_0 \oplus {\rosp(p\,|q)}_1 \; $  
by setting, for all  $ \, s \in \{0,1\} \, $,
$$ 
\; {\rosp(p\,|q)}_s  \, := \,  \Big\{ \ell \in \rgl(p\,|q) \;\Big|\; 
\phi\big(\ell(x),y\big) = -{(-1)}^{s\,|x|} \, \phi\big(x, \ell(y)\big) \; 
\forall \, x, y \in k^{p|q} \,\Big\}  
$$
and we call  $ \rosp(p\,|q) $  the {\it orthosymplectic\/}  Lie superalgebra.
Note that  $ \, \rosp(0|\,q) \, $  is the symplectic Lie algebra 
$\frak{sp}(q)$, 
while  $ \rosp(p|\,0) $  is the  orthogonal Lie algebra $\frak {so}(p)$.
                                                     \par
Again, all the  $ \rosp(p|\,q)$'s  are  Lie superalgebras 
{\sl of classical type\/}.  
Moreover, if  $ \, m, n \geq 2 \, $,  we have:
$$  
\displaylines{
   {\rosp(2m\!+\!1|2n)}_0 = \, \frak{so}(2m\!+\!1) \oplus \frak{sp}(2n) \; ,
\quad  {\rosp(2m|2n)}_0 = \, \frak{so}(2m) \oplus \frak{sp}(2n)  \cr
%
%
   {\rosp(p\,|2n)}_1 \, = \; f_p \otimes f_{2n}  \quad \forall \; p > 2 
\; ,  \;\;\qquad  {\rosp(2|2n)}_1 \, = \, f_{2n}^{\,\oplus 2}  }
$$

\end{examples}

We now introduce some terminology in order to be able to state
the classification theorem.

\begin{definition}
Define the following Lie superalgebras:
 \vskip4pt
\noindent
 \; {\it (1)} \;   $ A(m,n) := \rsl(m \! + \! 1 |\, n \! + \! 1) \, $,  \,
$ A(n,n) := \rsl(n \! + \! 1 |\, n \! + \! 1) \Big/ {k I_{2n}} \; $,   \,\;  $ \forall \; m \! \neq \! n \, $;
 \vskip4pt
\noindent
 \; {\it (2)} \hfill   $ B(m,n) := \rosp(2m+1|\,2n) \; $,   \hfill   $ \forall \;\;\; m \geq 0 \, $,  $ \, n \geq 1 \; $; \quad {\ }
  \vskip4pt
   {\it (3)} \hfill   $ C(n) := \rosp(2|\,2n-2) \, $,   \hfill   for all  $ \, n \geq 2 \, $; \quad {\ }
 \vskip4pt
   {\it (4)} \hfill   $ D(m,n) := \rosp(2m|\,2n) \; $,   \hfill   for all  $\, m \geq 2 \, $,  $ \, n \geq 1 \, $; \quad {\ }
 \vskip4pt
   {\it (5)} \hfill   $ \quad \displaystyle{
 P(n) := \, \left\{ \begin{pmatrix}
                         A  &  B  \\
                         C  & -A^t
\end{pmatrix} \in \rgl(n\!+\!1|\,n\!+\!1) \;\bigg|\;
  \begin{matrix}  A \in \rsl(n\!+\!1)  \\
           B^t = B \, , \; C^t = -C
  \end{matrix}   \;\right\} } $   \hfill {\ }
 \vskip4pt
   {\it (6)} \hfill   $ \displaystyle{
 Q(n) := \, \bigg\{\! \begin{pmatrix}
                         A  &  B  \\
                         B  &  A
\end{pmatrix} \in \rgl(n\!+\!1|\,n\!+\!1) \;\bigg|\;\, B \in \rsl(n\!+\!1) \;\bigg\} \Bigg/ k I_{2(n+1)} }. $  \hfill {\ }
\end{definition}

\begin{theorem} \label{classificationtheorem}
Let  $ k $  be an algebraically closed field of characteristic zero.  
Then the Lie superalgebras of classical type  are either isomorphic 
to a simple Lie algebra or to one of the following 
Lie superalgebras:
$$  
\begin{array}{c}
A(m,n) \, ,  \;\; m \! \geq \! n \! \geq \! 0 \, , \, m+n > 0 \, ;  \quad  
B(m,n) \, ,  \;\; m \geq 0, n \geq 1 \, ;   \quad  C(n) \, ,  \;\; n \geq 3 
\\ \\
D(m,n) \, ,  \;\; m \geq 2 , \, n \geq 1 \, ;  \qquad  
P(n) \, ,  \;\; n \geq 2 \, ;   \quad  Q(n) \, ,  \;\; n \geq 2  
\\ \\
F(4) \; ;  \qquad  G(3) \; ;  \qquad  
D(2,1;a) \, ,  \;\; a \in k \setminus \{0, -1\}. 
\end{array} 
$$
For the definition of $F(4)$, $G(3)$,    
$D(2,1;a)$, and for the proof, 
we refer to \cite{ka}.

\end{theorem}

\section{Chevalley basis 
and Kostant integral form}
\label{chevalleybasis}

The main ingredient to construct a Chevalley supergroup starting
from a complex Lie superalgebra $\fg$ of classical type
is the {\sl Chevalley basis}.
This is an homogeneous basis for $\fg$, consisting of elements
that have brackets expressed as {\sl integral} combinations
of the basis elements. Consequently a Chevalley basis determines what is
called the {\it Chevalley Lie algebra} $\fg_\Z$ of $\fg$, which
is an integral form of $\fg$.

\medskip

Assume $\fg$ to be a Lie superalgebra of classical type different
from $A(1,1)$, $P(3)$, $Q(n)$.
We want to leave out these pathological cases for which our construction
holds, but with a more complicated set of statements and proofs.
We invite the reader to go to \cite{fg} for a complete and
unified treatment of all of these cases. We also consider $D(2,1;a)$
for only integral values for the coefficient $a$.

\medskip

Let us fix a Cartan subalgebra $\fh$ of $\fg$, that is a maximal
solvable Lie subalgebra of $\fg$. The adjoint action of $\fh$ on
$\fg$ gives the usual {\sl root space} decomposition of $\fg$:
$$
\fg  \, = \, \fh \oplus \bigoplus_{\alpha \in \Delta} \fg_\alpha 
$$
where
$$
\fg_\alpha := \big\{\, x \in \fg \;\big|\; 
[h,x] = \alpha(h) x \, , \; \forall \; h \! \in \! \fh \big\}
$$
and $\Delta=\Delta_0 \cup \Delta_1$ with
$$
\begin{array}{c}
\Delta_0  \, := \,  \big\{\, \alpha \in \fh^* \setminus \{0\} \;\big|\; 
\fg_\alpha \cap \fg_0 \not= \{0\} \big\}  \, = \,  
\{\,\hbox{{\sl even roots\/}  of  $ \fg $} \,\}.  \\ \\
\Delta_1  \, := \,  \big\{\, \alpha \in \fh^* \;\big|\; 
\fg_\alpha \cap \fg_1 \not= \{0\} \big\}  \, = \,  
\{\,\hbox{{\sl odd roots\/}  of  $ \fg $} \,\}.  
\end{array}
$$
As in the ordinary setting we shall call $\Delta$ \textit{root system}
and the $\fg_\alpha$'s the \textit{ root spaces}.
If we fix a simple system (see \cite{ka} for its definition)
the root system splits into positive and negative roots, exactly as
in the ordinary setting:
$$
\Delta=\Delta^+ \coprod \Delta^-, \qquad
\Delta_0=\Delta_0^+ \coprod \Delta_0^-, \qquad
\Delta_1=\Delta_1^+ \coprod \Delta_1^-.
$$
\medskip

\begin{observation}

\begin{enumerate}

\item Notice that the definition allows 
$\Delta_0 \cap \Delta_1 \neq \emptyset$,
as in fact happens for $\fg=Q(n)$, where the roots are simoultaneously all even
and odd and the root spaces have all dimension $1|1$.

\medskip

\item $\Delta_0 $  is the root system of the reductive Lie algebra  
$ \fg_0 \, $,  while  $ \Delta_1 $  is the set of weights of the 
representation of  $ \fg_0 $  in  $ \fg_1 \, $.

\end{enumerate}

\end{observation}

\medskip

   If  $ \fg $  is not of type  $P(n)$  or  $Q(n)$,  
there is an even non-degenerate, invariant bilinear form on  
$\fg$,  whose restriction to  $ \fh $  is in turn an invariant 
bilinear form on  $ \fh \, $.  On the other hand, if $ \fg $  
is of type  $ P(n) $  or  $ Q(n) $,  then such a form on  
$ \fh $  exists because  $ \fg_0 $  
is simple (of type  $ A_n $), though it does not come by restricting
an invariant form on the whole $\fg$.
  
\medskip

If $ \big( x, y \big) \, $ denotes such form,  we can identify  $ \fh^* $  
with  $ \fh $,  via  $ H'_\alpha \mapsto \big(H'_\alpha, \big)$.
We can then transfer $\big( ,  \big)$ to $\fh^*$ in the natural way: 
$\big( \alpha, \beta \big) = \big( H'_{\alpha}, H'_{\beta} \big) \; $.
Define $H_\alpha:= 2{H'_\alpha \over \big(H'_\alpha, H'_\alpha\big)}$
when the denominator is non zero. 
When $\big(H'_\alpha, H'_\alpha\big)=0$ such renormalization can
be found in detail in \cite{ik}.
We call $H_\alpha$ the \textit{coroot} associated with $\alpha$.   

\medskip

We summarize in the next proposition all the relevant properties of
the root system, sending the reader to \cite{ka,sc,se} for the
complete story.

\begin{proposition}  \label{prop-root-syst}
Let  $ \fg $  be a Lie superalgebra of  {\sl classical type}, as above,
that is $\fg \neq A(1,1)$, $P(3)$, $Q(n)$,
and  let $ \, n \in \N \, $.

\medskip

\textit{(a)} \;  
$\; \Delta_0 \cap \Delta_1 = \emptyset \;$,  

\medskip

\textit{(b)} \;  $ -\Delta_0 = \Delta_0 \;$,  
$\; -\Delta_1 \subseteq \Delta_1 \; $.    
If  $ \, \fg \not= P (n) \; $,  then  $ \; -\Delta_1 = \Delta_1 \; $.

\medskip

\textit{(c)} \;  
Let  $ \, \fg \not= P(2) \, $,  and  
$ \, \alpha $,  $ \beta \in \Delta \, $,  
$ \; \alpha = c \, \beta \, $,  \, with  
$ \, c \in \KK \setminus \{0\} \, $.  Then
$$  
\alpha , \beta \in \Delta_r \quad (r=0,1)  
\, \Rightarrow \,  c = \pm 1 \;\; ,  \,\quad\,
\alpha \in \Delta_r \, , \, \beta \in \Delta_s \, ,  
\; r \not= s  \, \Rightarrow \,  c = \pm 2 \;\; .  
$$

\medskip

\textit{(d)} \;   
$ \; \dim_\KK(\fg_\alpha) = 1 \; $  for each  $ \, \alpha \in \Delta \, $.  
\end{proposition}

We are finally ready to give the definition of {\sl Chevalley basis}.

\begin{definition}  \label{def_che-bas}
We define a  {\it Chevalley basis\/}  of  a Lie superalgebra $ \fg $  
as above any homogeneous  basis  
$$
B = \big\{ H_1 \dots H_l, \, \, X_\alpha, \, \, \alpha \in \Delta \big\}
$$  
of $\fg$ as complex vector space, with the following requirements:

 \textit{(a)}  
$ \, \big\{ H_1 , \dots , H_\ell \big\} \, $  is a  
basis of the complex vector space $\fh$.  
Moreover
$$
\fh_\Z  \, := \,  \text{\it Span}_{\,\Z} \big\{ H_1 , \dots , H_\ell \big\}  
\, = \,  \text{\it Span}_{\,\Z} \big\{ H_\alpha \,\big|\, \alpha \! 
\in \! \Delta \}.
$$
%
 \vskip8pt
%
 \textit{(b)}  \hskip4pt   $ \big[ H_i \, , H_j \big] = 0 \, ,   \hskip9pt
 \big[ H_i \, , X_\alpha \big] = \, \alpha(H_i) \, X_\alpha \, ,   \hskip15pt  \forall \; i, j \! \in \! \{ 1, \dots, \ell \,\} \, ,  \; \alpha \! \in \! \Delta \; $;
 \vskip11pt
   \textit{(c)}  \hskip7pt   $ \big[ X_\alpha \, , \, X_{-\alpha} \big]  \, = \,  \sigma_\alpha \, H_\alpha  \hskip25pt  \forall \;\; \alpha \in \Delta \cap (-\Delta) $
 \vskip4pt
\noindent
 with  $ H_\alpha $  as after  \ref{prop-root-syst},  and  $ \; \sigma_\alpha := -1 \; $  if  $ \, \alpha \in \Delta_1^- \, $,  $ \; \sigma_\alpha := 1 \; $  otherwise;
 \vskip13pt
   \textit{(d)}  \quad  $ \, \big[ X_\alpha \, , \, X_\beta \big]  \, = \, c_{\alpha,\beta} \; X_{\alpha + \beta}  \hskip17pt   \forall \;\, \alpha , \beta \in \Delta \, : \, \alpha \not= -\beta$,  \, with $c_{\al,\be}  \in \Z$.
More precisely,

\begin{itemize}
\item If $(\al,\al) \neq 0$, or $(\beta,\beta) \neq 0$,
$c_{\al,\be} =\pm(r+1)$ or (only if $\fg=P(n)$), $c_{\al,\be} =\pm(r+2)$,
where $r$ is the length of the $\alpha$-string through $\beta$.

\item If  $(\al,\al) \neq 0=(\beta,\beta)=0$, $c_{\al,\be} =\be(\al)$.
\end{itemize}
\end{definition}

Notice that this definition clearly extends to direct sums of finitely 
many  of the $ \fg $'s under the above hypotheses.

\begin{definition}
If  $ B $  is a Chevalley basis of  a Lie superalgebra $ \fg \, $ as above,  
we set  
$$
\fg_\Z := \text{span}_\Z\{B\} \subset \fg
$$ 
and we call it the  {\it Chevalley superalgebra\/}  of  $ \fg $.
\end{definition}

Observe that $ \fg_\Z $  is a Lie superalgebra over  $ \Z \, $ inside $\fg$.
Since a Chevalley basis  $ B $  of $\fg$ is unique up to a choice of a sign 
for each root vector and the choice of the  $ H_i $'s
we have that  $ \fg_\Z $  is independent of the choice of  $ B \, $
(but of course depends on the choice of $\fh$ the Cartan subalgebra). 

\medskip

The existence of a Chevalley basis for the families  $ A $,  $ B $,  
$ C $,  $ D $  is a  known result; for example an almost
explicit Chevalley basis for types  $ B \, $,  $ C $  and  $ D \, $ 
is in \cite{sw}, while for $A$ is a straightforward calculation.  
More in general, an abstract existence result, 
with a uniform proof, is given in  \cite{ik}  for all  {\sl basic\/}  types.
In \cite{fg} we provide an existence theorem for {\sl all} cases
giving both a case by case analysis, comprehending all
Lie superalgebras of classical type and a uniform proof, 
that however leaves out the $P(n)$ case.

\medskip

We now turn to another important ingredient for our construction:
the {\sl Kostant $\Z$-form}.

\medskip

\begin{definition}  \label{def-kost-superalgebra}
Let $\fg$ be a complex Lie superalgebra of classical type over $\C$ and
let $B=\big\{ H_1 \dots H_\ell, X_\alpha, \alpha \in \Delta \big\}$
be a Chevalley basis.
We define  the \textit{ Kostant superalgebra},  $ \kzg \, $,
the  $ \Z $--superalgebra  inside  $ U(\fg) \, $, generated by
$$  
X_\alpha^{(n)} \; ,  \!\quad X_\gamma \; ,  
\!\quad  {\textstyle \Big(\! {H_i \atop n} \!\Big)}
\qquad \hfill \forall \;\; \alpha \! \in \! \Delta_0 \; , 
\; n \! \in \! \N \, ,
\; \gamma \! \in \! \Delta_1 \, , \; i = 1, \dots, \ell \, ,
$$
where
$$
X_\alpha^{(n)} \! := X_\alpha^n \big/ n! \qquad   
\bigg(\! {H \atop n} \!\bigg)  \; :=  \;  
{\frac{\, H (H \! - \! 1) \cdots (H \! - \! n \! + \! 1) \,}{n!}}
\in U(\fg)  
$$
for all $H$ in $\fh$. These are called respectively \textit{divided
  powers} and \textit{binomial coefficients}.
\end{definition}

Notice that we can remove all the binomial coefficients
corresponding to coroots $H_i$'s relative
to even roots and still generate the superalgebra $\kzg$.
In fact a classical result (see \cite{st} pg 9) tells us that
the even divided powers generate all such binomial coefficients.
Unfortunately we cannot obtain the odd coroot binomial coefficients
and this is because the $X_\ga$, for $\ga \in \Delta_1$ appear
only in degree one.

\medskip

As in the ordinary setting (see \cite{st} pg 7) we have a 
PBW type of result for $ \kzg \, $ providing us with a
$\Z$-basis for the Kostant superalgebra. The proof is very
similar to the ordinary setting and we send the reader to
\cite{fg} for more details.

\begin{theorem}  \label{PBW-Kost}
 The Kostant superalgebra  $ \, \kzg \, $  is a free\/  $ \Z $--module.  
For any given total order  $ \, \preceq \, $  of the set  
$ \, {\Delta} \cup \big\{ 1, \dots, \ell \big\} \, $,  
a  $ \, \Z $--basis  of  $ \kzg $  is the set   
of ordered ``PBW-like monomials'', i.e.~all products without repetitions:
of factors of type:  
$$ 
X_\alpha^{(n_\alpha)}, \quad  
\Big(\! {H_i \atop n_i} \!\Big), \quad 
X_\gamma 
$$
$ \, \alpha \in \Delta_0 \, $,  $ \, i \in \big\{ 1, \dots, \ell \big\} \, $,  
$ \, \gamma \in \Delta_1 \, $  and  $ \, n_\alpha $,  
$ n_i  \in \N $  ---   taken in the right order with respect to  
$ \, \preceq \; $.
\end{theorem}

\section{Chevalley supergroups} \label{chesupergroups}

This section is devoted to the construction
of {\sl Chevalley supergroups} and to prove they are
supergroup schemes. 

\medskip

Let  $ \fg $ be a complex Lie superalgebra of classical type,  
$B=\big\{ H_1 \dots H_\ell$, $X_\alpha$, $\alpha \in \Delta \big\}$ 
a Chevalley basis of  $ \fg $  and $ \kzg $ its 
Kostant superalgebra.
We start with a finite dimensional
complex representation $V$ of $\fg$ and
the notion of {\sl admissible lattice} in $V$.

\begin{definition} \label{rational}
Let  $V$  be a complex finite dimensional 
representation for $ \fg $. We say that $V$ is  {\it rational\/}  if  
$\fh_\Z := \text{\it Span}_\Z\big(H_1,\dots,H_\ell\big)$  
acts diagonally on  $ V $  with integral eigenvalues.

\medskip

Notice that this condition is automatic for semisimple
Lie algebras, while it is actually restrictive for some
Lie superalgebras as the next example shows.

\medskip

\begin{example}
Let $\fg=\rsl(m|n)$ and $\fh$ the diagonal matrices, so that 
$\fh_\Z= \text{\it Span}_\Z\big\{E_{m,m}+E_{m+1,m+1},\, 
E_{ii}-E_{i+1,i+1}, i \neq m \big\}$, where $E_{ij}$ denotes
an elementary matrix. Let $V$ be a representation
with highest weight $\Lambda=\lambda_1 \ep_1+ \dots +\lambda_m \ep_m+
\mu_1 \delta_1 + \dots \mu_n \delta_n$, where
$\ep_i: \fh \lra \C$, $\ep_i(E_{jj})=\delta_{ij}$ and
similarly for $\delta_k$.
We have that (see \cite{ka}) $V$ is finite dimensional if
and only if $\lambda_i-\lambda_{i+1}$, $\mu_j-\mu_{j+1} \in \Z^+$,
$i=1 \dots m-1$, $j=1 \dots n-1$,
in other words if and only if $\Lambda(H_i) \in \Z^+$ for $i \neq m$.
There are hence no conditions on $\Lambda(H_m)=\lambda_m+\mu_1$.
Consequently if we pick any (non integral) complex number for such
a sum and we build the induced module, we shall obtain a finite
dimensional representation for $\fg$ where $H_m$ acts diagonally,
with a complex, non integral eigenvalue.
\end{example}

Let us now fix $V$ a finite dimensional
rational complex semisimple representation of $\fg$.

\medskip

We say that an integral lattice  $ M $ in $V$ is 
{\it admissible} if it is  $ \kzg $--stable.
\end{definition}

As in the ordinary setting any rational complex finite dimensional
semisimple representation of $\fg$ admits an admissible lattice $M$, which 
is generated by the highest weight vector $v$ 
and it is the sum of its weight components $M_\mu$. In particular
if $V$ is simple we have:
$$
M=\kzg \cdot v, \qquad M=\oplus M_\mu.
$$

The next proposition establishes the existence of an integral
form of $\fg$ stabilizing the admissible lattice $M$ inside
the representation $V$. We send the reader to \cite{fg} \S 5
for the proof.

\begin{theorem}  \label{stabilizer}
Let  $\fg$, $V$ and $M$ as above. Define:
$$
\fg_V = \big\{ X \! \in \! \fg \,\big|\, X.M \subseteq M \big\}.
$$ 
If  $ \, V $  is faithful, then
$$  
\fg_V  \, = \,  \fh_V \; 
{\textstyle \bigoplus \big( \oplus_{\alpha \in \Delta}} \,
\Z \, X_{\alpha} \big) \;\; ,  \qquad  \fh_V \, := \,
\big\{ H \in \fh \;\big|\; \mu(H) \in \Z \, ,
\,\; \forall \; \mu \in \Lambda \big\}  
$$
where  $ \Lambda $  is the set of all weights of  $ \, V $.  
In particular,  $ \fg_V $  is a lattice in  $ \fg \, $,  
and it is independent of the choice of the admissible lattice  
$ M $  (but not of course of $ \, V $).
\end{theorem}

We end this discussion by saying that $\fg_\Z$ corresponds to
the adjoint representation of $\fg$ and that in general all the
integral forms $\fg_V$ lie between the two integral forms
$\fg_{roots}$ and $\fg_{weights}$
corresponding respectively to the root and the fundamental
weight representations:
$$
\fg_{roots} \subset \fg_V \subset \fg_{weights}.
$$

We now start the construction of the Chevalley supergroup associated
with the data $\fg$ and $V$.

\medskip

Let $k$ be a generic field.

\medskip

Definition \ref{rational} and Theorem \ref{stabilizer} 
allow us to move from the
complex field to a generic field quite easily as the next definition
shows.

\begin{definition}
Let $\fg$ be a complex
Lie superalgebra of classical type 
(as usual $\fg \neq A(1,1)$, $P(2)$, $Q(n)$).
Let $V$ be a faithful rational 
complex representation of $\fg$, $M$ an admissible lattice
in $V$.

\medskip

Define:
$$
\fg_k :=  k \otimes_\Z \fg_V, \qquad  
V_k := k \otimes_\Z M, \qquad  
U_k(\fg) := k \otimes_\Z \kzg.
$$
\end{definition}

We are now ready to the define the super equivalent of the
one-parameter subgroups in the classical theory.
As we shall see, homogeneous one-parameter subgroups appear in the 
super setting with
three different dimensions: $1|0$, $0|1$ and $1|1$. In order to keep
the analogy with the ordinary setting, we neverthless have preferred
to keep the terminology {\sl one-parameter subgroup}, though in the
supersetting the term ``one'' can be misleading.

\begin{definition} \label{one-par}
Let $X_\al$, $X_\be$, $X_\ga$ be root vectors in the Chevalley basis,
$\al \in \Delta_0$, $\be$, $\ga$ $\in \Delta_1$, with
$[X_\be,X_\be]=0$, $[X_\ga,X_\ga] \neq 0$.
 
We define \textit{homogeneous one-parameter subgroups}
the following supergroup functors from the categories
of superalgebras to the category of sets:
$$  
\begin{array}{rl}
x_\alpha(A) \!  &  := \big\{ \exp\!\big( t \, X_\alpha \big) \;\big|
\; t \in A_0 \,\big\} \; = \\ \\
& =\big\{ \big( 1 + t \, X_\alpha + 
t^2 \, {X_\alpha^2 \over 2}+ \cdots \big) \;\big| \; t \in A_0 \,\big\} 
\subset \rGL(V_k)(A), 
\\ \\ 
 x_\beta(A) \!  &  :=   
\big\{ \exp\!\big( \vartheta \, X_\beta \big)
\;\big|\; \vartheta \in A_1 \,\big\}  \; 
=\big\{ \big( 1 +  \vartheta \, X_\beta \big) 
\;\big|\; \vartheta \in A_1 \,\big\} \subset  \rGL(V_k)(A), 
\\  \\
x_\gamma(A) \!  &  := \,  
\big\{ \exp\!\big( \vartheta \, X_\gamma + 
t \, X_\gamma^{\,2} \big) \;\big|\; \vartheta \in A_1 \, , 
\, t \in A_0 \,\big\}  \, =  \\ \\
\phantom{\bigg|}  &  \phantom{:}= 
\,  \big\{ \big( 1+\vartheta \, X_\gamma\big) 
\exp\!\big( t \, X_\gamma^{\,2} \big) \;\big|\; 
\vartheta \in A_1 \, , \, t \in A_0 \,\big\}\subset  \rGL(V_k)(A). 
\end{array}  
$$
\end{definition}

Notice that the infinite sums reduce to finite ones since $X_\al$ and
$X_\be$ act as nilpotent operators on $V_k$. As usual we identify 
a generic root vector $X_\al$ with its image under the representation
of $U_k(\fg)$ in $V_k$ (the divided powers come at hand exactly at this
point).

\medskip

One can readily see that the functors $x_\al$, $x_\be$, $x_\ga$ 
are representable, hence
they are algebraic supergroups in the sense of \S \ref{supergeometry} and their
representing Hopf superalgebras are respectively 
$k[x]$, $k[\xi]$ and $k[x,\xi]$, where as usual the roman
letters correspond to even elements while the greek letters to odd ones.
The comultiplication is coadditive except for the element 
$x$ in $k[x,\xi]$: 
$x \mapsto 1 \otimes x + x \otimes 1 + \xi \otimes \xi$.
It is very clear by looking at the Hopf
superalgebras representing $x_\al$, $x_\be$ and $x_\ga$  that
the superdimensions of these supergroups are respectively $1|0$,
$0|1$ and $1|1$. It is not hard to see that these are all of the allowed
superdimensions for homogeneous one-parameter subgroups (see \cite{fg}
for more details). 
\medskip

As an abuse of notation we shall sometimes write for $t \in A_0$
and $\theta \in A_1$:
$$
\begin{array}{rl}
x_\al(t) &:=\exp\!\big( t \, X_\alpha \big), \\ \\ 
 x_\beta(\theta) &:=   \exp\!\big( \vartheta \, X_\beta \big) 
= 1 +  \vartheta \, X_\beta  \\ \\ 
x_\gamma(\bt) &:=   \exp\!\big( \vartheta \, X_\gamma + 
t \, X_\gamma^{\,2} \big), \qquad \bt=(t,\theta).
\end{array}
$$

\medskip

We now turn to the definition of the generators of what classically
is the {\sl maximal torus}.

\begin{definition}
For any  $ \, \alpha \in \Delta \subseteq \fh^* \, $,  let  
$ \, H_\alpha \in \fh_\Z \, $  as in \ref{prop-root-syst}.  
Let  $ \, V = \oplus_\mu V_\mu \, $  be the splitting of  $ V $  
into weight spaces. As  $ V $  is rational, we have  
$ \, \mu(H_\alpha) \in \Z \, $  for all  $ \, \alpha \in \Delta \, $.
Define:  
$$ 
\; h_\alpha(t).v  \, := \,  t^{\mu(H_\alpha)} \, v  \in V_k(A)\; 
\, \forall v \in (V_k)_\mu \, ,   \, \mu \in \fh^* \, \quad t \in A^\times,
\quad A \in \salg
$$
Notice that this defines an operator $h_\al(t) \in \rGL(V_k)(A)$. 
Hence we can define:
$$
h_H(t):=
\prod h_{\alpha}^{a_\alpha}(t) \in \rGL(V_k)(A), \qquad 
H=\sum a_\alpha H_\alpha.
$$
\end{definition}

We have immediately that $h_H$ defines a supergroup functor: 
$$
h_H:\salg \lra \sets, \qquad 
h_H(A):=\left\{h_H(t) \quad \big| \quad t \in A^\times\right\}
$$
which is clearly representable, its Hopf superalgebra being
given by $k[x,x^{-1}]$ with comultiplication $x \mapsto x \otimes x$.

\medskip

We now want to define an ordinary algebraic group associated
with the ordinary Lie algebra $\fg_0$, the even part of $\fg$.
One must exert some care at this point, since a Chevalley basis
for $\fg_0$ is not in general the even part or even a subset 
of a Chevalley basis for $\fg$,
even if $\fh=\fh_0$,
that is the Cartan subalgebras for $\fg_0$ and $\fg$ coincide.
Let us illustrate the first of these phenomenons with simple example.
Let us look at $A(2,1)$. 
$\fh_\Z=(\fh_\Z)_0={\mathrm Span_\Z} \{H_1, H_2, H_3, H_4\}$,
with $H_i=E_{ii}-E_{i+1,i+1}$, $i \neq 3$, $H_3=E_{33}+E_{44}$,
where the $E_{ij}$'s are the elementary matrices with $1$ in the
$(i,j)$-th position and zero elsewhere. We have only one odd coroot $H_3$.
This is an even vector, 
that we would however miss if we were to consider just the
even coroots, that is the coroots corresponding to even root spaces.
These are the coroots of the roots in $\Delta_0$ the root system
associated with $A_2 \oplus A_1$ the {\sl simple} even part 
of $\fg$, which in this case is not
the same as $\fg_0=A_2 \oplus A_1 \oplus \C$, which is {\sl reductive}.
To produce an instance of the second phenomenon is more complicated.
The point is that we can have that the span of the odd coroots
may contain some of the even coroots and consequently we can omit
those even coroots, so that the Chevalley basis will not be a subset
of the Chevalley basis of the even part. This happens for
example in the $D(m,n)$ case.
 
\medskip

It is possible to construct a reductive algebraic
group $\bG_0$ overcoming these difficulties.
$\bG_0$ will encode also the information contained in the extra odd
coroot (it is in fact possible always to reduce to the case of just one
odd coroot) and such construction is explained in detail in \cite{fg}.
The group $\bG_0$ is constructed following Chevalley's 
phylosophy, but taking into account the extra odd coroot, which
would be otherwise missing.
On local superalgebras $\bG_0$ is described as follows. Let $\bG_0'$ be
the ordinary algebraic group scheme associated with the semisimple
part of $\fg$ (which could be smaller than $\fg_0$ as in the $A(m,n)$
case) and let $T:\salg \lra \sets$, $T(A)=\langle h_H(A) \quad | \quad
H \in \fh_\Z\rangle$. This is in general larger than the maximal torus
$T_0$ in $G_0$, since it contains the extra odd root (though one must
be aware of some exceptions as we detail in the observation below). 
If $A$ is a local superalgebra we define:
$$
\bG_0(A)= \langle \bG_0'(A), T(A)  \rangle.
$$
It is possible to show that this definition extends to any
superalgebra $A$ and that the functor so obtained is representable
(see \cite{fg} section 5).

\begin{observation} \label{osp} 
We want to observe that there are cases in which the missing odd
root can be somehow recovered without extra work. Let us look at
the example of $\rosp(1|2)$. The roots are $\al$, $2\al$ and the
corresponding coroots are $H_\al=2H_{2\al}$, $H_{2\al}$ (the relation 
between coroots depends on the chosen normalization). Consequently,
we have that by taking just the even coroot $H_{2\al}$, we can get
both the coroots $H_\al$ and $H_{2\al}$, 
so in this case it is not necessary to add anything
more, in other words $\bG_0=\bG_0'$. 
Clearly this phenomenon is observed for all the superalgebras
$B(m,n)$. Notice that the even coroot $2\al$ corresponds to
the adjoint representation of the even part $\rsl_2$ of $\rosp(1|2)$.
As we shall see in our construction, this will tell us that we cannot
obtain  from $\rosp(1|2)$ a Chevalley supergroup whose reduced group
is ${\mathrm{SL}}_2$ and in fact we shall see that there is only one algebraic supergroup
associated with $\rosp(1|2)$ and its reduced algebraic group is
${\mathrm {PSL}}_2$. 
This fact has consequences on the questions regarding
which supergroups can be built using our method and we plan to
fully explore this in a forthcoming paper.
\end{observation}

\medskip

We are finally ready for the definition of {\sl Chevalley
supergroup functor}.

\begin{definition} \label{def_Che-sgroup_funct} 
Let  $\fg$ be a complex Lie superalgebra of classical type
and  $V$ a faithful rational complex representation of $\fg$.  
We call  the {\it Chevalley supergroup},  associated to  $\fg$  
and  $ V \, $,  the functor  $ \; \bG : \salg \lra \text{(grps)} \; $ 
defined as:
$$  
\begin{array}{rl}
\bG(A) &:= \Big\langle\, \bG_0(A) \, , \, x_\beta(A) \,\;
\Big|\;\, \beta \in \Delta_1 \, \Big\rangle=  \\ \\
&=\Big\langle\, \bG_0(A) \, , \, 1+\theta_\beta X_\beta \,\;
\Big|\;\, \beta \in \Delta_1, \, \theta_\be \in A_1 \Big\rangle
\subset \rGL\big(V_k(A)\big).
\end{array}
$$
In other words $\bG(A)$  is the subgroup of  
$\rGL\big(V_k(A)\big)$  generated by  $ \bG_0(A) $  described above
and the $0|1$ one-parameter subgroups  $ x_\beta(A) $  with  
$ \, \beta \in \Delta_1 \, $.
$\bG$ is defined on the arrows in the natural way, since $\bG(A)$ is
a subgroup of $\rGL\big(V_k(A)\big)$.

\medskip

From the classical theory (see \cite{sga} 5.7) we know that on
local algebras, since $\bG_0$ is reductive:
$$
\bG_0(A)=\Big\langle\, x_\alpha(A), h_i(A) \, \big| \, 
i = 1, \dots, \ell \, , \; \alpha \in \Delta_0 \Big\rangle.  
$$
Consequently on local superalgebras we then have:
$$
\bG(A)=\Big\langle\, x_\alpha(A), h_i(A) \, \big| \, 
i = 1, \dots, \ell \, , \; \alpha \in \Delta \Big\rangle.  
$$
We call  {\it Chevalley supergroup functor\/} the functor 
$ \; G : \salg \lra \text{(grps)} \; $  defined as:
$$
G(A)=\Big\langle\, x_\alpha(A), h_i(A) \,  \big| \,
i = 1, \dots, \ell \, , \; \alpha \in \Delta \Big\rangle.  
$$
\end{definition}

In \cite{fg} we explore more deeply the relation between the two
functors $\bG$ and $G$ and we show that $\bG$ is the
{\sl sheafification} of $G$. This important property sheds
light on our construction and it is actually needed in the key proofs,
since it provides a more explicit way to handle the Chevalley
supergroups. Neverthless, given the scope of the present work, 
we shall not give the definition of sheafification of a functor,
in order to avoid the technicalities involved, that are not
adding any insight into our construction.
For all the details we send
the reader to the appendix in \cite{fg} and, for the ordinary setting to 
\cite{vst}, where the sheafification of functors is fully explained
and to \cite{sga} 5.7.6 for its application to reductive groups.

\medskip

The fact that we have defined the Chevalley supergroup $\bG$
as a functor does not automatically imply that  it  is  {\it representable},
in other words, that  it  is the functor of points of an
algebraic supergroup scheme. This is a new question specific to the
supersetting, in fact in the ordinary setting, the definition of Chevalley group
is given only on fields, the group is exhibited an abstract
group and only later one shows it is has an algebraic scheme structure.
On the other hand in the supergeometric environment looking
at superobjects on fields only will not give us much information since 
the odd coordinates disappear when we look at points over a
field, thus leaving us with just the underlying ordinary group.
In other words $\bG(k)=\bG_0(k)$ for all fields $k$, since
the $\theta_\be$'s in Definition \ref{def_Che-sgroup_funct}
are nilpotent.

\medskip

In order to prove the representability of $\bG$, we shall
give a series of lemmas regarding $G$, which is more accessible
than $\bG$, since we know its generators for all $A \in \salg$.
As in the ordinary setting the key to the theory are the explicit
formulas for the commutators. The proof is
a straightforward generalization of the corresponding proofs for
the ordinary setting (see \cite{st} \S 3), which we stated as $(1)$
of \ref{comm_1-pssg}.

\medskip

Before this 
in order to properly state our results and the intermediate steps to
obtain them, we need to define the following auxiliary sets.

\begin{definition}  \label{subgrps}
 For any  $ \, A \in \salg \, $,  we define the subsets of  $ G(A) $
  $$  \begin{array}{rl}
   G_1(A)  &  \! := \;  \left\{\, {\textstyle \prod_{\,i=1}^{\,n}} \,
x_{\gamma_i}(\vartheta_i) \;\Big|\; n \in \N \, , \; \gamma_1, \dots, \gamma_n \in \Delta_1 \, , \; \vartheta_1, \dots, \vartheta_n \in A_1 \,\right\}  \\
  \\
   G_0^\pm(A)  &  \! := \;  \left\{\, {\textstyle \prod_{\,i=1}^{\,n}}
\, x_{\alpha_i}(t_i) \;\Big|\; n \in \N \, , \; \alpha_1, \dots, \alpha_n \in \Delta^\pm_0 \, , \; t_1, \dots, t_n \in A_0 \,\right\}  \\
  \\
   G_1^\pm(A)  &  \! := \;  \left\{\, {\textstyle \prod_{\,i=1}^{\,n}} \,
x_{\gamma_i}(\vartheta_i) \;\Big|\; n \in \N \, , \; \gamma_1, \dots,
\gamma_n \in \Delta^\pm_1 \, , \; \vartheta_1, \dots, \vartheta_n \in A_1
\,\right\}  \\
   \\
   G^\pm(A)  &  \! := \;  \left\{\, {\textstyle \prod_{\,i=1}^{\,n}} \,
x_{\beta_i}(\bt_i) \;\Big|\; n \! \in \! \N \, , \, \beta_1, \dots, \beta_n \! \in \! \Delta^\pm , \, \bt_1, \dots, \bt_n \! \in \! A_0 \! \times \! A_1 \right\}
=  \\
   \phantom{\Bigg|}  &  \hfill   \, = \; \big\langle G_0^\pm(A) \, , \, G_1^\pm(A) \big\rangle
      \end{array}  $$
 \vskip-9pt
\noindent
 Moreover, fixing any total order  $ \, \preceq \, $  on $ {\wDelta}_1^\pm \, $,  and letting  $ \, N_\pm = \big| {\wDelta}_1^\pm \big| \, $,  we set
  $$  G_1^{\pm,<}(A)  \,\; := \;  \left\{\; {\textstyle \prod_{\,i=1}^{\,N_\pm}} \, x_{\gamma_i}(\vartheta_i) \,\;\Big|\;\, \gamma_1 \prec \cdots \prec \gamma_{N_\pm} \in \Delta^\pm_1 \, , \; \vartheta_1, \dots, \vartheta_{N_\pm} \in A_1 \,\right\}  $$
and for any total order  $ \, \preceq \, $  on  $ {\wDelta}_1 \, $,  and letting  $ \, N := \big| \Delta \big| = N_+ + N_- \, $,  we set
  $$  G_1^<(A)  \,\; := \;  \left\{\, {\textstyle \prod_{\,i=1}^N} \, x_{\gamma_i}(\vartheta_i) \,\;\Big|\;\, \gamma_1 \prec \cdots \prec \gamma_{\scriptscriptstyle N} \in \Delta_1 \, , \; \vartheta_1, \dots, \vartheta_N \in A_1 \,\right\}  $$
   \indent   Note that for special choices of the order, one has  $ \; G_1^<(A) = G_1^{-,<}(A) \cdot G_1^{+,<}(A) \; $  or  $ \; G_1^<(A) = G_1^{+,<}(A) \cdot G_1^{-,<}(A) \; $.
                                                                      \par
%
\end{definition}

\begin{remark}
 Note that  $ G_1(A) $,  $ G_0^\pm(A) $,  $ G_1^\pm(A) $  and  $
 G^\pm(A) $  are subgroups of  $ G(A) \, $,  while  $ G_1^{\pm,<}(A) $
 and  $ G_1^<(A) $  instead are  
{\sl not},  in general.  
%
\end{remark}

\begin{lemma}  \label{comm_1-pssg} 

\begin{enumerate}

\item Let  $ \; \alpha, \beta \in \Delta_0 \, $,  $ \, A \in \salg \, $  and  
$ t,u \in A_0 \, $.  
Then there exist  $ \, c_{ij} \! \in \! \Z \, $  such that
$$  
\big( x_\al(t) \, , \, x_\be(u) \big)  \, = \,  
{\textstyle \prod} \, x_{i \, \alpha + j \, \beta}\big(c_{ij} \, 
t^iu^j \big)  \;\; \in \;\;  G_0(A).  
$$


\item Let  $ \; \alpha \in \Delta_0 \, $,  
$ \, \gamma \in \Delta_1 \, $,  $ \, A \in \salg \, $  and  
$ \, t \in A_0 \, $,  $ \, \vartheta \in A_1 \, $.  
Then there exist  $ \, c_s \! \in \! \Z \, $  such that
  $$  \big( x_\gamma(\vartheta) \, , \, x_\alpha(t) \big)  \, = \,  
{\textstyle \prod_{s>0}} \, x_{\gamma + s \, \alpha}\big(c_s \, 
t^s \vartheta\big)  \;\; \in \;\;  G_1(A),
$$
(the product being finite).  
More precisely, with  $ \, \varepsilon_k = \pm 1 \, $  and  $ \, r \in
\Z\, $,
$$  
\big( 1 + \vartheta \, X_\gamma \, , \, x_\alpha(t) \big)  \, = \,
{\textstyle \prod_{s>0}} \, \Big( 1 + 
{\textstyle \prod_{k=1}^s \varepsilon_k \cdot {{s+r} \choose r}} 
\cdot t^s \vartheta \, X_{\gamma + s \, \alpha} \Big)  
$$
where the factors in the product are taken in any order (as they do commute).

%

\item Let  $ \; \gamma , \delta \! \in \! \Delta_1 \, $,  
$ \, A \! \in \! \salg \, $,  $ \, \vartheta, \eta \! \in \! A_1 \, $.  
Then (notation of  Definition \ref{def_che-bas})
$$  
\big( x_\gamma(\vartheta) \, , \, x_\delta(\eta) \big)  \; = \;\,  
x_{\gamma + \delta}\big(\! - \! c_{\gamma,\delta} \; \vartheta \, \eta
\big)  
\; = \;  \big(\, 1 \! - \! c_{\gamma,\delta} \; \vartheta \, \eta \,
X_{\gamma + \delta} \,\big)  
\;\; \in \;\;  G_0(A)  
$$
if  $ \; \delta \not= -\gamma \; $;  otherwise, for  $ \; \delta = -\gamma \, $,  we have
$$  
\big( x_\gamma(\vartheta) \, , \, x_{-\gamma}(\eta) \big)  \; = \;
\big(\, 1 \! - \vartheta \, 
\eta \, H_\gamma \,\big)  \; = \;  h_\gamma\big( 1 \! - \vartheta \,
\eta \big)  
\;\; \in \;\;  G_0(A).  
$$
%
\item Let  $ \; \alpha , \beta \in \Delta \, $,  $ \, A \in \salg \,
  $,  
$ \, t \in U(A_0) \, $,  $ \, \bu \in A_0 \! \times \! A_1 = A \, $.  
Then
$$  
\hskip31pt   h_\alpha(t) \; x_\beta(\bu) \; {h_\alpha(t)}^{-1}  \; =
\;  
x_\beta\big( t^{\beta(H_\alpha)} \, \bu \big)  \;\; \in \;\;
G_{p(\beta)}(A)  
$$
where $p(\beta)$ denotes as  usual the {\sl parity} of a root $\beta$,
that is $p(\be)=0$ if $\be \in \Delta_0$ and $p(\be)=1$ if $\be \in
\Delta_1$. 
\end{enumerate}
\end{lemma}

We are still under the simplifying assumption $\fg \neq Q(n)$ hence
$\Delta_0 \cap \Delta_1=\emptyset$. We stress that our results hold
for {\sl all} Lie superalgebras of classical type, but we choose in
the present work for clarity of exposition to restrict ourselves to
$\fg \neq A(1,1), P(3), Q(n)$.

\medskip

As a direct consequence of the commutation relations, we have
the following proposition involving the sets we have introduced:
$G^\pm$, etc. The proof is a simple exercise.

\begin{theorem} \label{g0g1}
Let  $ \, A \in \salg \, $.  There exist set-theoretic factorizations
$$  
\begin{array}{c}
G(A)  \; = \;  G_0(A) \; G_1(A) \; = \;  G_1(A) \; G_0(A)  \\ \\
G^\pm(A)  \; = \;  G^\pm_0(A) \; G^\pm_1(A) \; = \;  G^\pm_1(A) \; G^\pm_0(A).
\end{array}  
$$
\end{theorem}

This decomposition has a further refinement that we state down below,
whose proof is
harder and we send the reader to \cite{fg} \S 5.3 for the details.

\begin{theorem} \label{realcrucial}
For any  $ \, A \in \salg \, $  we have
$$  
\begin{array}{c}
G(A)  \; = \;  G_0(A) \, G_1^<(A) \; = \;  G_1^<(A) \, G_0(A)  
\end{array}
$$
\end{theorem}

From the previous results we have that a generic $g \in G(A)$ can
be factorized  (once we choose a suitable ordering on the roots):
$$
g=g_0g_1^+g_1^-, \qquad g_0 \in G_0(A), \quad g_1^\pm = G^{\pm,<}_1(A).
$$

\medskip

The next theorem gives us the key to the representability of $\bG$,
by stating the uniqueness of the above decomposition.
Again for the proof see \cite{fg}, 5.3.

\begin{theorem} \label{crucial}
Let the notation be as above.
For any  $ \, A \in \salg \, $,  the group product gives 
the following bijection:
$$  
G_0(A) \times G_1^{-,<}(A) \times G_1^{+,<}(A) 
\lhook\joinrel\relbar\joinrel\relbar\joinrel\relbar\joinrel\twoheadrightarrow
\, G(A)  
$$
and all the similar bijections obtained by permuting the factors  
$ G_1^{\pm,<}(A) $  and the factor  $ G_0(A) $.
\end{theorem}

As one can readily see, the functors $G_1^{\pm,<}: \salg \lra \sets$
are representable and they are the functor of points of an odd
dimensional affine super space:
$G_1^{\pm,<} \cong \mathbb A^{0|N^\pm}$, for 
$N^\pm=|\Delta_1^\pm|$. Then this, together with the 
definition of $\bG$ gives:
$$
\bG \cong \bG_0 \times G_1^{-,<} \times G_1^{+,<} =
\bG_0 \times \bA^{0|M}
$$
for $M=N^++N^-$.
Consequently $\bG$ is representable, since it is the direct product of
representable functors.
We have sketched the proof of the main result of the paper:

\begin{theorem} \label{rep-chev}
The Chevalley supergroup $\bG:\salg \lra \sets$,
$$
\bG(A) := \Big\langle\, \bG_0(A) \, , \, x_\beta(A) \,\;
\Big|\;\, \beta \in \Delta_1 \, \Big\rangle
$$
is representable.
\end{theorem}

\medskip

The next proposition establishes how much the Chevalley supergroup
scheme $\bG$ we have built depends on the chosen representation.
It turns out that two different complex $\fg$-representations $V$ and $V'$
(as in beginning of Sec. \ref{chesupergroups}),
with weight lattices $L_{V'} \subset L_{V}$
of the same complex Lie superalgebra $\fg$ of classical type give raise to
a morphism between the corresponding
Chevalley supergroups, with kernel inside the
center of $\bG$, as it happens in the ordinary setting. This is actually
expected, since the kernel is related with the fundamental group, which
is a topological invariant, unchanged by the supergeneralization.

\begin{theorem} \label{mainthm}
Let  $ \bG $  and  $ \bG' $  be two Chevalley supergroups constructed
using faithful complex representations $ V $  and  $ V' $  of the
same  complex Lie superalgebra of classical type $ \fg $.  
Let $ L_V $,  $ L_{V'} $  be the corresponding lattices of weights.             
If  $ \, L_V \supseteq L_{V'} \, $,  then there exists a unique morphism  $ \; \phi : \bG \longrightarrow \bG' \; $  such that  $ \,\; \phi_A\big(1 + \vartheta \, X_\alpha\big) = 1 + \vartheta \, X'_\alpha \; $,  \, and  $ \,\; \text{\it Ker}\,(\phi_A) \subseteq Z\big(\bG(A)\big) \, $,  \; for every local algebra  $ A \, $.  Moreover,  $ \phi $  is an isomorphism if and only if  $ \, L_V = L_{V'} \, $.
\end{theorem}

We observe that this theorem tells us that our construction
of $\bG$ does not depend on the chosen representation $V$, but only
on the weight lattice of $V$. 
In particular $\bG$ is independent of the choice of an admissible lattice.
\medskip

In the end we want to ask the following question: does our
construction provide all the algebraic supergroups whose Lie
superalgebra is of classical type? The answer to this question is
positive and we plan to explore furtherly the topics in a forthcoming
paper.

\section{Examples and further topics} \label{examples}

In this final section we want to discuss some examples and
to indicate possible further developments and applications
of the theory we have described.

\medskip

We start by discussing how our construction can be generalized to
other Lie superalgebras, provided some conditions are satisfied.

\medskip

We list down below some requirements a Lie superalgebra must satisfy
so that we can try to replicate our construction.

\medskip

We start from a complex Lie superalgebra 
$\fg=\langle X_a \, \big| \, a \in {\mathcal A} \rangle$ 
generated (as Lie superalgebra) by the homogeneous elements $X_a$,
where $a \in {\mathcal{A}}$ a finite set of indices,
and a complex finite dimensional representation $V$.

\medskip

We assume the following.

\begin{enumerate} 

\item $\fg$ admits a basis 
$B \supset \{X_a\}_{a \in {\mathcal{A}}}$ and an integral form 
$\fg_\Z=\mathrm{span}_\Z \{B\}$ in which 
all the brackets are integral combinations of elements in $B$;

\item There exists a suitable $\Z$-subalgebra of $U(\fg)$ denoted by
$\uzg$ 
$\subset U(\fg)$ 
admitting a PBW theorem. In other words, $\uzg$ is a free $\Z$-module
with a basis consisting of suitable monomials, which form also a
basis for $U(\fg)$.

\item $V$ contains an integral lattice $M$ stable under $\uzg$;

\item There is well defined algebraic group $\bG_0$ over $\Z$, whose 
$k$-points embed into  
$\rGL(k \otimes M)$ and whose corresponding Lie algebra is $\fg_0$.
This will allow us to consider its functor of points
$\bG_0:\salg \lra \sets$. 

\end{enumerate} 

If the requirements listed above are satisfied, then we can
certainly give the same definition as in \ref{def_Che-sgroup_funct}
and define the {\sl Chevalley-type supergroup functor}.

\medskip

Notice that the first part, up to section \ref{chevalleybasis} is devoted
to prove $(1)$--$(3)$ for $\fg$ of classical type. $(4)$ for the
classical type is discussed in
\cite{fg} section 4. 

\begin{definition}
Let $\fg$ and $V$ be as above. Define as before (compare before
\ref{one-par}):
$$
\fg_k :=  k \otimes \fg_V, \qquad  
V_k := k \otimes M, \qquad  
U_k(\fg) := k \otimes \uzg.
$$
Define the \textit{Chevalley-type supergroup functor} as the functor
$\bG: \salg \lra \sets$ given on the objects by:
$$
\bG(A)=\big\langle \bG_0(A), \quad 1+\theta X_b, 
\quad X_b \quad \hbox{odd}, \quad \theta \in A_0 \big\rangle
\subset \rGL(V_k)(A).
$$
In other words $\bG(A)$ is the subgroup of $\rGL(V_k)$ generated by
the $A$-points of $\bG_0(A)$ and the elements $1+\theta X_b$.
Again we identify an element $X_b$ with its image in the
representation $V_k$.
\end{definition}

As we have already remarked after Definition \ref{def_Che-sgroup_funct},
this definition does not ensure $\bG$ to be representable,
hence to be rightfully called a supergroup scheme,
and in fact a key role in the proof of the representability
of this functor in the case of $\fg$ of classical type, is played by the
commutation relations between the elements generating 
the group $G(A)$. 

\medskip

Before we go to the representability issues, let us give
an example of Lie superalgebra together with a class of representations,
which is {\sl not} of classical type
and yet it satisfies the requirements listed above, hence it
admits a Chevalley-type supergroup functor.

\begin{example}
Let us consider the Heisenberg Lie superalgebra $\fH$, which is
generated by an even generator $e$ and by $2n$ odd generators
$a_i$, $b_i$, $i=1 \dots n$, with the only non zero brackets:
$$
[a_i,b_j]=\delta_{ij} e, \qquad i,j=1 \dots n.
$$
Define the following irreducible
faithful complex 
representation (see \cite{ka} \S 1.1).
Let $V=\wedge(\xi_1 \dots \xi_n)$, the complex exterior algebra with
generators $\xi_1 \dots \xi_n$.
$V$ is a complex representation
for $\fH$ by setting:
$$
e \cdot u = \al u, \quad 
a_i \cdot u={{\partial u} \over {\partial \xi_i}}, \quad
b_i \cdot u=\al \xi_i u, \qquad  \al \in \C.
$$

\medskip

Assume we take $\al \in\Z$.

\medskip

If we set $\{X_a\}_{a=1, \dots,  2n}=\{a_i,b_i, \, i=1 \dots n\}$ and
$B=\{e,a_i,b_i, \, i=1 \dots n \}$ we have immediately satisfied item $(1)$. 
As for item $(2)$, we have that: 
$$
\{\left( \begin{array}{c} e \\ n \end{array} \right) a_{i_1} \dots a_{i_p}
b_{j_1} \dots b_{j_q}\}, \qquad 1 \leq i_1 < \dots < i_p \leq n, \quad
1 \leq j_1 < \dots < j_q
$$
is a $\Z$-basis of 
$$
U_\Z(\fH):=\big\langle \left( \begin{array}{c} e \\
        n \end{array} \right), a_i, b_j \big\rangle \subset U(\fH).
$$
$V$ contains the following integral lattice stable under $U_\Z(\fH)$:
$$
M={\mathrm span}_\Z \left\{\xi_{i_1} \dots \xi_{i_m} \quad
\big| \quad 1 \leq  i_1 < \dots < i_m \leq n \right\}  
$$

Finally certainly $a_i$ and $b_i$ act as nilpotent operators
and consequently also item $(3)$ is satisfied.

\medskip

As for item $(4)$, 
we have immediately
that $\bG_0 \cong k$ is an algebraic group, the additive group of
the affine line.
In the representation $V_k$ the elements in $\bG_0(k)$ act as follows:
$$
h_e(t) \cdot u=t^\alpha u
$$
hence $\bG_0(k)$ is embedded into $\rGL(V_k)$ as the diagonal matrices:
$$
\bG_0(k)=\left( \begin{array}{cccc} t^\al & 0 & \dots & 0 \\
\vdots & & & \vdots \\
0 & 0 & \dots & t^\al \end{array} \right) \subset \rGL(V_k).
$$
Its functor of points is hence given simply by taking $t \in A_0$:
$$
G_0(A)=\left( \begin{array}{cccc} t^\al & 0 & \dots & 0 \\
\vdots & & & \vdots \\
0 & 0 & \dots & t^\al \end{array} \right) \subset \rGL(V_k), \qquad t \in A_0.
$$

\end{example}

We now go to the problem of representability of the
functor $\bG$ we have defined. 
Besides the technical problems involved, our proof of
the representability issue for the Chevalley supergroup functor
discussed in \ref{chesupergroups} relies on the following facts:

\begin{enumerate}

\item $\bG(A)=\bG_0(A)G_1^<(A)$, where 
$G_1^<(A)=\{(1+\theta_{a_1}X_{a_1}) \dots (1+\theta_{a_n}X_{a_n})\}$
where $a_{i}<a_{i+1}$ in an order on $\cal A$, the index set of the
indices $a_i$.

\item The above decomposition is
unique, that is $\bG(A)= \bG_0(A) \times G_1^<(A)$ is unique.

\item $G_1^< \cong \bA^{0|N}$ for a suitable $N$.

\end{enumerate}

Clearly this leads immediately to the representability of
the functor $\bG$, since it is the direct product of two
representable functors.

\medskip

Coming back to our example of the Heisenberg superalgebra, 
by a direct calculation very similar to the one in section
\ref{chesupergroups} one sees that the commutator:
$$
(1+\theta X_a, 1+\eta X_b)=1+c \theta \eta e, \quad c \in \Z.
$$ 
Notice that $1+c \theta \eta e$ acts on $u \in V_k$ as
a diagonal matrix with entries $1+c\theta \eta t^\al$. This is
an element in $G_0(A)$ since 
$(1+c \theta \eta t)^\al=1+c \theta \eta t^\al$.
By repeating the reordering arguments as in \cite{fg} 5.15 and 5.16 
one can show that properties (1)--(3) are
satisfied, hence giving us the representability of the
Chevalley-type supergroup functor for
the Heisenberg Lie superalgebra. Consequently we have define the
Heisenberg supergroup associated to the Heisenberg Lie superalgebra
in the following way:
$$
\bG(A)=\langle \bG_0(A), 1+\theta_i a_i, 1+\eta_i b_i \rangle
\subset \rGL(V)(A).
$$ 

\medskip

We now want consider an important question: What are the algebraic
supergroups, 
that we can construct using this method, in other words using
the Chevalley supergroup construction? 
One could be tempted to say that we can obtain always an algebraic supergroup
whose reduced group corresponds to a reductive group
with Lie algebra the even part of the given Lie superalgebra of classical type.
The situation, however, is more complicated
and this is not always the case. Let us examine an interesting
example.

\medskip

Using the method of Super Harish-Chandra pairs (SHCP) (see for
example \cite{koszul}) we have that there exists a simply 
connected\footnote{A \textit{simply connected} supergroup is
a supergroup with simply connected underlying topological space.} 
Lie supergroup corresponding to the Lie superalgebra
$\rosp(p|q)$. Since $\rosp(p|q)=\rso(p) \otimes \rsp(q)$, this
supergroup can be rightfully called the \textit{spin supergroup}
and it is very important in physics.
For example in \cite{vsv} ch. 5, Varadarajan discusses its construction
for arbitrary signatures.

\medskip

However, with our construction, it is not always possible to reach
the simply connected supergroup associated with the given Lie superalgebra.
Let us clarify this phenomenon with an example.

\begin{example}
Let us consider: 
$$
\rosp(1|2):=\left\{
\begin{pmatrix} 0 & \al \quad \be \\
\begin{array}{c} \be \\ -\al \end{array}& B
\end{pmatrix} | B \in \rsl_2\right\} \subset \rgl(1|2).
$$
A basis is given by:
$$
\begin{array}{c}
e=\begin{pmatrix} 0 & 0 & 0 \\
0 & 0 & 1 \\ 0 & 0 & 0 \end{pmatrix}, 
f=\begin{pmatrix} 0 & 0 & 0 \\
0 & 0 & 0 \\ 0 & 1 & 0 \end{pmatrix}, 
h=\begin{pmatrix} 0 & 0 & 0 \\
0 & 1 & 0 \\ 0 & 0 & -1 \end{pmatrix}, \\ \\ 
x=\begin{pmatrix} 0 & 0 & 1 \\
1 & 0 & 0 \\ 0 & 0 & 0 \end{pmatrix},
y=\begin{pmatrix} 0 & 1 & 0 \\
0 & 0 & 0 \\ -1 & 0 & 0 \end{pmatrix} 
\end{array}
$$
where ${\mathrm span}_\C
\{e,f,h\} \cong \rsl_2$ and the other non zero brackets are:
$$
\begin{array}{cccc}
[e, y] = -x, &   [f,x] = -y, & [x,x]=2e, & [y,y]=-2f, \\
\end{array}
$$
$$
\begin{array}{cccc}
[x,y]=h,  &  [h,x]=x, & [h,y]=-y. &
\end{array}
$$
It is clear that a Cartan subalgebra over $\C$ can be chosen as
$\fh={\mathrm span}_\C \{h\}$ and a Chevalley basis for $\rosp$ is
$$
{\mathcal B}=\{x,y,e,f,h\}.
$$

Now we go to the choice of the module $V$, that determines
the Chevalley supergroup. In order to obtain a supergroup with
simply connected underlying topological space, we must choose
for $\rsl(2) \subset \rosp(1|2)$ the fundamental representation,
that is $\fg_V=\{h/2,e,f\}$. However as one can see right away,
such choice is forbidden, since otherwise we would have:
$$
[h/2,x]=x/2, \qquad [h/2,y]=-y/2.
$$
Consequently the only Chevalley supergroup associated with
$\rosp(1|2)$ has as underlying topological space the
adjoint group for $\mathrm{SL}(2)$, that is $\mathrm{PSL}(2)$.

\medskip

This shows that our method will not yield the spin supergroup, and
we believe that this Lie supergroup is not algebraic. We plan to
explore this and further topics in this section in a forthcoming
paper.
\end{example}


\begin{thebibliography}{99}

\bibitem{am} M.~F.~Atiyah, I.~G.~MacDonald,  
{\it Introduction to Commutative Algebra},  Addison-Wesley 
Publ.~Comp.~Inc., London, 1969.

\bibitem{bkl} J.~Brundan, A.~Kleshchev,  
{\it Modular representations of the supergroup  $ Q(n) $,  I},  
J.~Algebra  {\bf 206}  (2003), 64--98.

\bibitem{bku} J.~Brundan, J.~Kujava,  
{\it A New Proof of the Mullineux Conjecture},  
J.~Alg.~Combinatorics  {\bf 18}  (2003), 13--39.

\bibitem{borel} A.~Borel,  {\it Properties and linear representations 
of Chevalley groups},  in:  {\sl Seminar on Algebraic Groups and Related 
Finite Groups},  A.~Borel et al.~(eds.), Lecture Notes in  Math.~{\bf 131},  
Springer, Berlin, 1970, pp. 1--55.


\bibitem{cf} C.~Carmeli, L.~Caston, R.~Fioresi, 
{\it  Mathematical Foundation of Supersymmetry}, 
with an appendix with I. Dimitrov, EMS Ser. Lect. Math., European
Math. Soc., Zurich, 2011.

\bibitem{dm} P.~Deligne, J.~Morgan,  {\it Notes on supersymmetry 
(following J.~Bernstein)},  in: ``Quantum fields and strings. 
A course for mathematicians'', Vol.~1, AMS, 1999.

\bibitem{dg} M.~Demazure, P.~Gabriel,  {\it Groupes Alg\'ebriques, Tome 1},  
Mason$\&$Cie \'editeur, North-Holland Publishing Company, 
The Netherlands, 1970.

\bibitem{eh} D.~Eisenbud, J.~Harris,  {\it The Geometry of Schemes},  
Graduate Texts in  Math.~{\bf 197},  Springer-Verlag, New York-Heidelberg, 
2000.

\bibitem{fg} R. Fioresi, F. Gavarini, {\it Chevalley Supergroups},
preprint arXiv:0808.0785, to be published on  Memoirs of the AMS, 2008.

 \bibitem{sga} M.~Demazure, A.~Grothendieck, {\it Sch{\'e}mas en groupes, III},
 S{\'e}minaire de G{\'e}om{\'e}trie Alg{\'e}brique du Bois Marie, Vol. 3, 1964.

\bibitem{fss} L.~Frappat, P.~Sorba, A.~Sciarrino,  {\it Dictionary on 
Lie algebras and superalgebras},  Academic Press, Inc., San Diego, CA, 2000.


\bibitem{ha} R.~Hartshorne,  {\it Algebraic geometry},  Graduate
Texts in Math.~{\bf 52}, Springer-Verlag, New York-Heidelberg,
1977.

\bibitem{hu} J.~E.~Humphreys,  
{\it Introduction to Lie Algebras and Representation Theory},  
Graduate Texts in Math.~{\bf 9},  Springer-Verlag, New York, 1972.

\bibitem{ik} K.~Iohara, Y.~Koga,  
{\it Central extensions of Lie Superalgebras},
Comment.~Math.~Helv.~{\bf 76}  (2001), 110--154.

\bibitem{ja} J.~C.~Jantzen,  {\it Lectures on Quantum Groups},  
Grad.~Stud.~Math.~{\bf 6},  Amer.~Math.~Soc., Providence, RI, 1996.

\bibitem{ka} V.~G.~Kac,  {\it Lie superalgebras},  
Adv.~Math.~{\bf 26}  (1977), 8--26.


\bibitem{koszul}  J.-L., Koszul, 
\textit{Graded manifolds and graded {L}ie algebras},
 {Proceedings of the international meeting on geometry and physics 
({F}lorence, 1982)}, 71--84, Pitagora,
Bologna, 1982.


\bibitem{ma} Y.~Manin,  {\it Gauge field theory and complex geometry},  
Springer-Verlag, Berlin, 1988.

\bibitem{ms} A.~Masuoka,  {\it The fundamental correspondences in 
super affine groups and super formal groups},  J.~Pure Appl.~Algebra  
{\bf 202}  (2005), 284--312.

\bibitem{sc} M.~Scheunert,  {\it The Theory of Lie Superalgebras},  
Lecture Notes in Math.~{\bf 716},  Springer-Verlag, 
Berlin-Heidelberg-New York, 1979.

\bibitem{se} V.~Serganova,  {\it On generalizations of root systems}, 
Comm.~Algebra  {\bf 24} (1996), 4281--4299.

\bibitem{st} R.~Steinberg,  {\it Lectures on Chevalley groups},  
Yale University, New Haven, Conn., 1968.

\bibitem{sw} B.~Shu, W.~Wang,  
{\it Modular representations of the ortho-symplectic supergroups},  
Proc.~Lond.~Math.~Soc.~(3)  {\bf 96}  (2008), no.~1, 251--271.

\bibitem{vst} A.~Vistoli,  {\it Grothendieck topologies, fibered 
categories and descent theory}, in: Fundamental algebraic geometry,  
1--104, Math.~Surveys Monogr.  {\bf 123},  Amer.~Math.~Soc., 
Providence, RI, 2005.

\bibitem{vsv} V.~S.~Varadarajan,  {\it Supersymmetry for
    mathematicians: an introduction},  Courant Lecture Notes  {\bf 1},
  AMS, 2004.

\end{thebibliography}
\end{document}